\documentclass[12pt,a4paper]{article}
\usepackage[utf8x]{inputenc}
\usepackage[english]{babel}
\usepackage[T1]{fontenc}
\usepackage{amsmath}
\usepackage{amsfonts}
\usepackage{amssymb}
\usepackage{makeidx}
\usepackage{graphicx}
\usepackage[left=3cm,right=2cm,top=3cm,bottom=2cm]{geometry}
\usepackage{setspace}
\usepackage{amsthm}
\usepackage{cite}

\usepackage[colorlinks=true,urlcolor=blue,
citecolor=red,linkcolor=blue,linktocpage,pdfpagelabels,
bookmarksnumbered,bookmarksopen]{hyperref}

\addto\captionsenglish{}

\newtheorem{definition}{Definition}[section]

\newtheorem{proposition}{Proposition}[section]

\newtheorem{theorem}{Theorem}[section]

\newtheorem{example}{Example}[section]

\newtheorem{lemma}{Lemma}[section]

\newtheorem{observation}{Remark}[section]

\newtheorem{corollary}{Corolary}[section]

\numberwithin{equation}{section}

\newcommand{\bo}{\begin{observation}}
\newcommand{\eo}{\end{observation}}
\newcommand{\bd}{\begin{definition}}
\newcommand{\ed}{\end{definition}}
\newcommand{\bp}{\begin{proposition}}
\newcommand{\ep}{\end{proposition}}
\newcommand{\bt}{\begin{theorem}}
\newcommand{\et}{\end{theorem}}
\newcommand{\bc}{\begin{corollary}}
\newcommand{\ec}{\end{corollary}}
\newcommand{\bl}{\begin{lemma}}
\newcommand{\el}{\end{lemma}}
\newcommand{\be}{\begin{example}}
\newcommand{\ee}{\end{example}}
\newcommand{\beq}{\begin{equation}}
\newcommand{\eeq}{\end{equation}}
\newcommand{\beqa}{\begin{equation*}}
\newcommand{\eeqa}{\end{equation*}}
\newcommand{\R}{\mathbb{R}}
\newcommand{\RN}{\mathbb{R}^{N}}
\newcommand{\N}{\mathbb{N}}
\newcommand{\Rdois}{\mathbb{R}^{2}}

\newcommand{\Ldois}{ L^{2}(\mathbb{R}^2) }

\newcommand{\Linf}{ L^{\infty}(\mathbb{R}) }

\newcommand{\B}{\mathcal{B}}
\newcommand{\I}{\mathcal{I}}
\newcommand{\A}{\mathcal{A}}
\newcommand{\F}{\mathcal{F}}
\newcommand{\K}{\mathcal{K}}

\newcommand{\IC}{\hookrightarrow}
\newcommand{\CF}{\rightharpoonup}

\newcommand{\Hum}{H^{1}(\mathbb{R}^2)}
\newcommand{\Ls}{ L^{s}(\mathbb{R}^2)}

\newcommand{\intR}{\displaystyle\int\limits_{\mathbb{R}^2}}
\newcommand{\un}{u_{n}}

\newcommand{\vn}{v_{n}}

\newcommand{\yn}{y_{n}}

\newcommand{\until}{\tilde{u}_n}
\newcommand{\RA}{\rightarrow}

\newcommand{\ds}{\displaystyle\int\limits}

\usepackage{times, color, xcolor}
\addto\captionsportuguese{
}

\begin{document}

 	\title{Existence and multiplicity results for a class of Kirchhoff-Choquard equations with a generalized sign-changing potential
		}
\author{
		Eduardo  de S. Böer$^{1}$ \footnote{ E-mail address: eduardoboer04@gmail.com, Tel. +55.51.993673377}, Ol\'{\i}mpio H. Miyagaki$^{1}$  \thanks{E-mail address: ohmiyagaki@gmail.com, Tel.: +55.16.33519178 (UFSCar)} \,\, and \,\, Patrizia Pucci$^{2}$ \thanks{Corresponding author} \footnote{E-mail address: patrizia.pucci@unipg.it}  \\
		{\footnotesize $^{1}$  Department of Mathematics, Federal University of S\~ao Carlos,}\\ {\footnotesize 13565-905 S\~ao Carlos, SP - Brazil}\\
		{\footnotesize $^{2}$  Dipartimento di Matematica e Informatica,
Universit\`{a} degli Studi di Perugia,}\\
		{\footnotesize 06123
Perugia, Italy}
		}
\noindent
				
	\maketitle

\begin{center}
{\emph{Dedicated to the memory of Professor Antonio Ambrosetti,
with high feelings of admiration\\for his notable contributions in Mathematics and great
affection}}
\end{center}

\noindent \textbf{Abstract:} In the present work we are concerned with the following Kirchhoff-Choquard-type equation $$-M(||\nabla u||_{2}^{2})\Delta u +Q(x)u + \mu(V(|\cdot|)\ast u^2)u = f(u)  \mbox{ in \ } \Rdois , $$ for $M: \R \RA \R \mbox{ \ given by \ } M(t)=a+bt$, $ \mu >0 $, $ V $ a sign-changing and possible unbounded potential, $ Q $ a continuous external potential and a nonlinearity $f$ with exponential critical growth. We prove existence and multiplicity of solutions in the \textit{nondegenerate} case and guarantee the existence of solutions in the \textit{degenerate} case.

\vspace{0.5 cm}

			
\section{Introduction}			
			
The present work is devoted to study existence and multiplicity of solutions to the following class of Kirchhoff-Choquard equations
\begin{equation}\label{P}
-M(||\nabla u||_{2}^{2})\Delta u +Q(x)u + \mu(V(|\cdot|)\ast u^2)u = f(u)\quad \mbox{ in \ } \Rdois ,
\end{equation}
where $ \mu > 0 $, $M: \R \RA \R$ is a Kirchhoff type function, $ Q: \Rdois \RA \R $ is a nonnegative potential, $V:\R \RA \R$ is a continuous sign-changing and possible unbounded potential and $ f:\R \RA \R $ is a continuous function with primitive $ F(t)=\int\limits_{0}^{t} f(s)ds $.

This paper was motivated by recent works dealing with Choquard equations with logarithmic kernel, such as \cite{[6], [10], [boer], [boer2]}, and some works of Kirchhoff-type equations, as for example \cite{[pucci], [olimpio]}. In the following, we make a quick literature overview.

On one hand, the following Choquard or nonlinear Schrödinger-Newton equation
\beq\label{i2}
-\Delta u + V(x) u + \gamma ( \Gamma_N \ast |u|^2) u = b|u|^{p-2}u,\ p> 2, \ b>0, \textrm{ \  \ in \ } \RN,
\eeq
where
$ \Gamma_N $ is the well-known fundamental solution of the Laplacian
$$
\Gamma_N (x) = \begin{cases} 
\dfrac{1}{N(2-N)\sigma_N}|x|^{2-N} &\textrm{ if \ } N\geq 3 ,\\
\dfrac{1}{2\pi}\ln |x|& \textrm{ if \ } N= 2 ,
\end{cases} 
$$
has been extensively studied in the case $ N=3 $, due to its relevance in physics. Although what the equation  name suggests, it was first studied by Fröhlich and Pekar in \cite{[12] , [11] , [22]}, to describe the quantum mechanics of a polaron at rest, in the particular case $ V(x)\equiv Constant > 0 $ and $ \gamma > 0 $. Then, in 1976, Choquard introduced the same equation in the study of an electron trapped in its hole. Moreover, Penrose has derived equation (\ref{i2}) while discussing about the self gravitational collapse of a quantum-mechanical system in  \cite{[18]}. See also \cite{[16]}.

In the case $ N=2 $, in \cite{[6]}, the authors have proved the existence of a ground state solution, using the Nehari  manifold and the existence of infinitely many geometrically distinct solutions, when $ Q: \Rdois \RA (0, \infty) $ is continuous and $ \mathbb{Z}^2 $-periodic, $ \mu > 0 $ and a nonlinearity of the form $ f(u)=b|u|^{p-2}u $, with $ b\geq 0 $ and $ p\geq 4 $. Then, intending to fill the gap, i.e., the situation $ 2<p<4 $, the paper \cite{[10]} deals with equation (\ref{P}) when $ Q(x)\equiv Constant > 0 $ and $ \mu >0 $, and $ f(u)=|u|^{p-2}u $, with $ 2<p<4 $, and provides  existence of a mountain pass solution as well as of a ground state solution. Finally, in \cite{[boer]}, the authors prove  existence and multiplicity results for the $p-$fractional Laplacian operator and in \cite{[boer2]}  existence and multiplicity results are derived for $ (p, N)$-Laplacian equations. Moreover, in \cite{[boer2]} the authors prove for the first time that, up to subsequence, Cerami sequences are bounded in the solution space.

Let us recall that, from a physical point of view, the local nonlinear terms on the right side of equation (\ref{i2}), such as $ b|u|^{p-2}u $, for $ b\in \R $ and $ p>2 $, usually appears in the Schrödinger equations as a way of modelling the interaction among particles.  We refer the reader to \cite{[boer]} for a complete overview in this topic.

On the other hand, the literature of Kirchhoff-type equations and its related elliptic problems is very interesting and quite large. As an example, we cite \cite{[jin]} where the authors consider the following equation
$$
\left\{\begin{array}{ll}
-\left(a+ b \ds_{\mathbb{R}^N}|\nabla u|^2 dx \right)\Delta u + u = f(x, u) \mbox{ in } \mathbb{R}^N ,\\[2ex]
u\in H^1(\mathbb{R}^N),
\end{array}
\right.
$$
and prove the existence of a sequence of radial solutions $(u_k)\subset H^1(\mathbb{R}^N)$ satisfying $ I(u_k)\RA   \infty $, as $ k \RA   \infty $. For a more detailed overview in the numerous results involving Kirchhoff equations, we refer the reader to \cite{[pucci], [olimpio], [pucci2], [liang]} and the references therein.

To finish, we emphasize  that nonlinearities with exponential behaviour appear frequently in applied problems, from physics to biology, which show us the importance of the studies on this topic. In this sense, we cite some works that deal with nonlinearities of Moser-Trudinger type \cite{[5], [Lam], [17], [olimpio]} and the references therein.

We intend to extend or complement the above mentioned works, considering a generalized sign-changing convolution potential, the Kirchhoff operator and a nonlinearity with critical exponential growth.

In the sequel we present the features of equation \eqref{P}. Throughout this paper, $ \R^+ = \{t\in \R \ ; \ t > 0\} $. In our work we are going to consider the following Kirchhoff function

\noindent $ (M) \
M: \R \RA \R \mbox{ \ given by \ } M(t)=a+bt \mbox{ , for all \ } t \in \R \mbox{ , with \ }a> 0 \mbox{ \ and \ } b \geq 0 \mbox{ or } a = 0 \mbox{ \ and \ } b >0.
$
The case where $ a>0 $ is called \textit{nondegenerate} while the situation in which $ a=0 $ is said to be \textit{degenerate}. We are going to consider both cases here.

Since our intention is to provide a way to solve problems with sign-changing potentials that can be unbounded from below, we require that $ V $ has a nontrivial negative part, $ V^- = \max\{-V, 0\} $. But some of the arguments can be modified in order to apply these techniques to positive potentials as well. The positive part of $V$ is defined as $\max\{ V, 0\}$. Thus, we assume that $V: \R^+ \RA \R$ is a real function verifying the following properties

\noindent $
(V_1) \ \mbox{There are real functions } a_1, a_2: \R^+ \RA \R \mbox{ such that } a_2 \in L^{\infty}(\R^+), a_{1, 0}= \inf\limits_{t\geq 2} a_1(t) > 0, a_{2, 0}=\inf\limits_{t \in \R^+}a_2(t) > 0 \mbox{ and }
$
$$
a_1(t)\ln(1+ t) \leq V^{+}(t) \leq a_2(t) \ln(1+t) , \forall \ t> 0.
$$

\noindent $(V_2)$ There exists a real function $ a_3: \R^+ \RA \R $ such that $a_3(t) > 0$ in a subset of $\R^+$ with positive measure,
$$V^{-}(t) \leq \dfrac{a_3(t)}{t}\quad \forall \ t> 0\quad\mbox{and}\quad
\begin{cases}
a_3 \in L^{\infty}(\R),\\
\mbox{or}\\
a_3(t)= t^{-\lambda}, \mbox{ for some } \lambda \in [1, 3) \mbox{ and  for all } t > 0,
\end{cases}
$$
$$
\mbox{ \ There exists an open subset } \I \subset \R^+ \mbox{ such that } V(t) < 0 \mbox{ for all } t \in \I.\leqno{(V_3)}
$$

Natural examples of potentials $ V $, satisfying conditions $ (V_1)-(V_3) $, have the following geometry: $ V(|x|)\RA  \infty $ as $ |x|\RA  \infty $ and, either its negative part is bounded or $ V(|x|)\RA - \infty $ as $ |x|\RA 0 $. Such behaviour is studied in Lemma \ref{l7}. Below, we present some examples of possible $ V $. Note that, in example \textbf{(b)}, the negative part is bounded.

\be\label{ex1}
\textbf{(a)} {\rm The most important example as a potential $ V $ is the logarithmic kernel, $ V(|x|)=\ln|x| $, for all $ x \in \Rdois $. Observe that condition $ (V_1) $ is satisfied with $ a_2(t)\equiv 1 $ and
$$
a_1(t) = \left\{ \begin{array}{lll}
\dfrac{\ln t}{2\ln(1+t)} \ \ , \mbox{ if } t \geq 2 , \\[2ex]
\dfrac{\ln 2}{2 \ln 3}(t-1) , \mbox{ if } 1 \leq t \leq 2 ,\\
0 \ \ \ \ \ \ \ \ \ \ \ \ \ \ \ \ \ \ , \mbox{ if } 0 \leq t \leq 1.
\end{array}
\right.  
$$
Conditions $(V_2)$ and $ (V_3) $ are verified with functions $ a_3, a_4, a_5 \equiv 1 $.

\noindent \textbf{(b)} Define $ V(|x|)=|x|^\alpha - |x|^\beta $, for $ 0< \beta < \alpha < 1 $ suitably chosen.

\noindent \textbf{(c)} One can also consider more exotic potentials, such as $ V:\R^+ \RA \R $ given by
$$
V(|x|)=
\begin{cases} 
-\dfrac{1}{|x|},&\mbox{if \ } 0<|x|\leq 1\\
2|x| - 3 ,&\mbox{if \ } 1\leq |x| \leq \dfrac{3}{2}\\
\ln\left(|x|-\dfrac{1}{2}\right),&\mbox{ if \ } \dfrac{3}{2}\leq |x|.
\end{cases}
$$}
\ee

\noindent
As for the external potential $ Q:\Rdois \RA \R $, we ask 
condition\smallskip

\noindent $ (Q) \
Q \in C(\Rdois , \R), \inf\limits_{x\in \Rdois}Q(x) = Q_0 > 0 \mbox{ and there exists } p \in (1,  \infty] \mbox{ such that } Q \in L^{p}(\Rdois) .
$

Finally, we recall that a function $ h $ has \textit{subcritical} exponential growth at $ +\infty $, if
$$
\lim\limits_{s\RA + \infty}\dfrac{h(s)}{e^{\alpha s^2}-1} = 0 \textrm{ \ , for all \ } \alpha >0 ,
$$
and we say that $ h $ has $ \alpha_0 $-\textit{critical} exponential growth at $ \infty $, if
$$
\lim\limits_{s\RA + \infty}\dfrac{h(s)}{e^{\alpha s^2}-1} = \left\{ \begin{array}{ll}
0, \ \ \ \forall \ \alpha > \alpha_{0} \\
\infty , \ \ \ \forall \ \alpha < \alpha_0
\end{array} \right. .
$$ 
Thus, inspired by works such as \cite{[ruf],[olimpio], [boer]}, we consider the following conditions over $ f $.

$$f\in C(\R , \R), f(0)=0 \mbox{ and has critical exponential growth with } \alpha_0 = 4\pi . \leqno{(f_1)}$$
$$ \lim\limits_{|t|\RA 0} \dfrac{|f(t)|}{|t|^\tau }=0 \mbox{, for some } \tau > 1.  \leqno{(f_2)}$$
$$ \mbox{ There exists }\ \theta \geq 4 \ \mbox{ such that}\   f(t)t \geq \theta F(t) > 0, \ \mbox{for all }\  t\in \R \setminus \{0\}. \leqno{(f_3)}$$
$$ \mbox{There  exist}\  q>4 \ \mbox{ and} \  C_q > 0 \ \mbox{ such that}\  F(t) \geq C_q |t|^q , \ \mbox{for all} \  t \in \R . \leqno{(f_4)}$$

\noindent
From conditions $ (f_1) $ and $(f_2)$, given $ \varepsilon >0 $, $\alpha > 4\pi,$ fixed, for all  $ p>2 $, we can find two constants $ K_1=K_1(p, \alpha , \varepsilon) > 0 $ and $ K_2=K_2(p, \alpha , \varepsilon) > 0$ such that
\begin{equation}\label{eq2}
f(t)\leq \varepsilon |t|^\tau +K_1 |t|^{p-1}(e^{\alpha t^2}-1) \ , \ \ \ \forall \ t \in \R,
\end{equation}
and
\begin{equation}\label{eq3}
F(t) \leq \varepsilon|t|^{\tau + 1} + K_2 |t|^p(e^{\alpha t^2}-1) \ , \ \ \ \forall t \in \R .
\end{equation}

\be\label{ex2}
{\rm As a prototype for nonlinearity $ f $ satisfying conditions $ (f_1)-(f_4) $, we can consider $ f: \R \RA \R $ given by 
$$
f(t) = C_q \left\{ \begin{array}{ll}
t^q \ \ \ \ \ \ \ \ \ \ \ \ \ , \mbox{ if } 0 \leq t \leq 1 \\
t^q e^{4\pi (t^2 -1)} , \mbox{ if } t > 1
\end{array} \right. ,
$$
for $ C_q > 0 $ sufficiently large and $ q> 3 $ and consider the its odd extension.}
\ee

\be\label{ex4}
{\rm From Example \ref{ex1}, one can see that problem \eqref{P} includes, as a very important particular case, the planar Schrödinger-Poisson system
$$
-\Delta u + Q(x) u + \mu (\ln|\cdot|\ast |u|^{2})u = f(u) \textrm{ \ in \ } \mathbb{R}^2.
$$}
\ee

We are now ready to enunciate our first main result.

\bt\label{t1}
Suppose $ (V_1)-(V_3) $, $ (Q) $, $ (f_1)-(f_4) $, $ a>0 $, $ b\geq 0 $, $ \mu > 0 $, $ q>4 $ and $ C_q>0 $ sufficiently large. Then,
\begin{itemize}
\item[(a)] problem \eqref{P} has a nontrivial solution at the mountain pass level, that is, there exists $ u\in X\setminus\{0\} $ such that $ u $ is a critical point for $ I $ and $ I(u)=c_{mp} $, where
\begin{equation}\label{eq16}
c_{mp} = \inf\limits_{\gamma \in \Gamma}\max\limits_{t \in [0, 1]}I(\gamma(t)),
\end{equation}
with $ \Gamma = \{\gamma\in C([0, 1], X) \ ; \ \gamma(0)=0 \mbox{ and } I(\gamma(1))< 0 \} $.

\item[(b)] Problem \eqref{P} has a nontrivial ground state solution, in the sense that, there is $ u\in X\setminus \{0\} $ that is a critical point to $ I $ and satisfies
$$
I(u)=c_g=\inf \{ I(v) \ ; \ v \in \K \}, \mbox{ \ where \ } \K = \{ v \in X \setminus\{0\} \ ; \ I'(v)=0\}.
$$
\end{itemize}
\et

Then, in order to get multiple solutions for \eqref{P}, we are going to apply a symmetric version of mountain pass theorem. To do so, we need to change condition $ (f_1) $ by the following.

$$f\in C(\R , \R), f(0)=0, f \mbox{ is odd and has critical exponential growth with } \alpha_0 = 4\pi . \leqno{(f_1')}$$

As a prototype example for this case, one can consider the odd extension of $ f $ given in Example \ref{ex2}.

\bt\label{t2}
Suppose $ (V_1)-(V_3) $, $ (Q) $, $ (f_1') $, $ (f_2)-(f_4) $, $ a>0 $, $ b\geq 0 $, $ \mu > 0 $, $ q>4 $ and $ C_q>0 $ sufficiently large. Then, problem \eqref{P} has infinitely many solutions.
\et

In the \textit{degenerate} case we also need some changes in the hypotheses for $ f $. First of all, in order to get the mountain pass geometry, we ask
$$ \lim\limits_{|t|\RA 0} \dfrac{|f(t)|}{|t|^\tau }=0 \mbox{, for some } \tau > 3.  \leqno{(f_2')}$$
Moreover, to obtain boundedness for Cerami sequences in $ \Hum $, we need
$$ \mbox{ There exists }\ \theta \geq 8 \ \mbox{ such that}\   f(t)t \geq \theta F(t) > 0, \ \mbox{for all }\  t\in \R \setminus \{0\}. \leqno{(f_3')}$$

\bt\label{t3}
Suppose $ (V_1)-(V_3) $, $ (Q) $, $ (f_1), (f_2'), (f_3'), (f_4) $, $ a=0 $, $ b> 0 $, $ q>4 $ and $ C_q>0 $ sufficiently large. Then,
\begin{itemize}
\item[(a)] there exists a value $ \mu_\ast >0 $ such that, for all $ \mu \in (0, \mu_\ast) $, problem \eqref{P} has a nontrivial solution at the mountain pass level, i.e., there exits $ u\in X\setminus\{0\} $ a critical point for $ I $ satisfying $ I(u)=c_{mp} $, where
\begin{equation*}
c_{mp} = \inf\limits_{\gamma \in \Gamma}\max\limits_{t \in [0, 1]}I(\gamma(t)),
\end{equation*}
with $ \Gamma = \{\gamma\in C([0, 1], X) \ ; \ \gamma(0)=0 \mbox{ and } I(\gamma(1))< 0 \} $.

\item[(b)] There exists a value $ \mu_{\ast \ast}\in (0, \mu_\ast]  $ such that, for all $ \mu \in (0, \mu_{\ast \ast}) $, problem \eqref{P} has a nontrivial ground state solution, in the sense that, there is $ u\in X\setminus \{0\} $ that is a critical point to $ I $ and satisfies
$$
I(u)=c_g=\inf \{ I(v) \ ; \ v \in \K \}, \mbox{ \ where \ } \K = \{ v \in X \setminus\{0\} \ ; \ I'(v)=0\}.
$$
\end{itemize}
\et

Throughout the paper, we consider the following notations: $ \Ls $ denotes the usual Lebesgue space with norm $ ||\cdot ||_s $; \ $ X' $ stands as the dual space of $ X $; \ $ B_r(x) $ is the ball centred in $ x $ with radius $ r>0 $, simply $ B_r $ if $ x=0 $;  \ $ r_1, r_2 $ will be real values verifying $ r_1, r_2 >1 $, $ r_1 \sim 1 $ and $ \frac{1}{r_1}+\frac{1}{r_2}=1 $; \ $ (y\ast u)(x) = u(x-y) $, for all $ x, y\in \Rdois $; \ $ x_n \searrow x $ will mean that $ x_n \RA x $ and $ x_n \geq x $, for all $ n\in \N $; \ $ K_i $, $ i\in \N $, denote  important constants present in the estimates; \  $ C_i $, $ i\in\N $, denote different positive constants whose exact values are not essential to the exposition of the arguments.

The paper is organized as follows: in Section 2 we present some technical results concerning the framework and boundedness of sequences in the solution space. Section 3 is devoted to analyse the geometry of the functional and the involved potentials. Moreover, we verify some boundedness and convergence properties. In Section 4 we consider the \textit{nondegenerate} case and present the proof of a key proposition and of Theorems \ref{t1} and \ref{t2}. Finally, in Section 5, we study the \textit{degenerate} case and prove the existence result related to it.

\section{Framework and Technical Results}

In this section, we are going to present  the space where we are going to look for solutions to \eqref{P} and some technical results concerning sequences in such space. To begin with, since we are going to use a variational approach, we introduce the Euler-Lagrange functional $I:\Hum \RA \R \cup \{\infty\}$ associated to \eqref{P},  given by
\begin{equation}\label{eq4}
I(u) = \dfrac{a}{2}\intR |\nabla u|^2 dx + \dfrac{b}{4}\left(\intR |\nabla u|^2 dx \right)^2 + \dfrac{1}{2}\intR Q(x) u^2(x) dx + \dfrac{\mu}{4}P(u) - \intR F(u) dx ,
\end{equation}
where $ P: \Hum \RA \R \cup \{\infty\} $ is defined as
\begin{equation}\label{eq5}
P(u) = \intR \intR V(|x-y|)u^2(x)u^2(y) dx dy .
\end{equation}
We also consider two auxiliary bilinear symmetric and positive forms $ \overline{P}_1 , \overline{P}_2 : \Hum \RA \R \cup \{\infty\} $ given by
\begin{small}
\begin{equation}\label{eq6}
\overline{P}_1(u) = \intR \intR V^+ (|x-y|)u(x)v(y) dx dy \mbox{ \ and \ } \overline{P}_2(u) = \intR \intR V^- (|x-y|)u(x)v(y) dx dy ,
\end{equation}
\end{small}
respectively, and the functionals $ P_1, P_2: \Hum \RA \R \cup \{\infty\} $ defined as $ P_1 = \overline{P}_1(u^2, u^2)$ and $ P_2 = \overline{P}_2(u^2 , u^2) $. Observe that $ P(u)=P_1(u)-P_2(u) $. 

Now, based on \cite{[6]}, we consider the slightly smaller Hilbert space
$$
X = \{ u\in \Hum \ ; \ ||u||_\ast < \infty \},\quad \mbox{where \ \ \ } ||u||_{\ast}^{2}= \intR \ln(1+|x|) u^2(x) dx ,
$$
endowed with the norm $ ||\cdot||_{X}^{2}=||\cdot||^2 + ||\cdot||_{\ast}^{2} $, where $ ||\cdot|| $ is the usual norm in $ \Hum $, and $ ||\cdot||_\ast $ comes from the inner product
$$
\langle u, v \rangle_\ast = \intR \ln(1+|x|) u(x) v(x)\ dx.
$$
Clearly $ ||\cdot|| \leq ||\cdot||_X $, so that $ X \IC \Hum \IC \Ls $, for all $ s \geq 2 $.

Our first aim in this section is to prove that $ I\in C^1(X, \R) $. Thus, we present some preliminary results to achieve this aim.

We start recalling the reader the well-known Moser-Trudinger inequality.

\bl\label{l4}
\cite{[5]} If $\alpha >0$ and $ u\in \Hum $, then
$$
\intR \left(e^{\alpha |u|^2} - 1 \right) \ dx  < \infty .
$$
Moreover, if $ ||\nabla u||_{2}^{2}\leq 1 $, $ ||u||_{2}^{2}\leq M < \infty $ and $ \alpha <  4 \pi $, then there exists a constant $ K_{\alpha, M}=K(M, \alpha) $, such that
$$
\intR \left(e^{\alpha |u|^2} - 1 \right) dx < K_{\alpha, M} .
$$
\el
Then, combined with \eqref{eq3} and Hölder inequality, for $ r_1, r_2 > 1 $, $ r_1 \sim 1 $ and $ \frac{1}{r_1}+\frac{1}{r_2}=1 $, we have
\begin{equation}\label{eq7}
\intR F(u) \ dx \leq \varepsilon ||u||_{\tau + 1}^{\tau + 1} + K_2 ||u||_{pr_2}^{p}\left(\intR (e^{r_1 \alpha|u|^2}-1) dx \right)^{\frac{1}{r_1}}   < 
\infty , \forall \ u \in \Hum .
\end{equation}
On the other side, from condition $ (V_1) $, there exists a constant $ K_3>0 $ such that
\begin{equation}\label{eq9}
P_1(u) \leq ||a_2||_\infty \intR \intR \ln(1+|x - y|) u^2(x)u^2(y) dxdy \leq K_3 ||u||_{\ast}^{2}||u||_{2}^{2},
\end{equation}
and, from condition $ (V_2) $ and the Hardy-Littlewood-Sobolev inequality (HLS) (found in \cite{[15]}), there is $ K_4 > 0 $ satisfying
$$
P_2(u) \leq ||a_3||_\infty \intR \intR \dfrac{1}{|x-y|}u^2(x)u^2(y) dx dy \leq K_4 ||u||_{\frac{8}{3}}^{4} ,
$$
if $ a_3 \in L^\infty (\R) $, and
\begin{equation}\label{eq10}
P_2(u) \leq \intR \intR \dfrac{1}{|x-y|^{\lambda+1}}u^2(x)u^2(y) dx dy \leq K_{HLS} ||u||_{\frac{8}{3-\lambda}}^{4},
\end{equation}
if $ a_3(t)=t^{-\lambda} $, for all $ t >0 $ and $ \lambda \in [-1, 3) $. One should observe that the constant $ K_4 $ also depend on the best Hardy-Littlewood-Sobolev constant, denoted by $ K_{HLS} $. From now on we consider only the second case, since   $a_3 \in \Linf $ can be treated similarly. Moreover, from condition $ (Q) $, we have, for $ p>1 $ given in $ (Q) $, with $ \frac{1}{p}+\frac{1}{p'}=1 $, that
\begin{equation}\label{eq13}
\intR Q(x) u^2 \ dx \leq ||Q||_p ||u||_{2p'}^{2} \mbox{ \ \ \ and \ \ \ } \intR Q(x)uv \ dx \leq ||Q||_p ||u||_{2p'}||v||_{2p'} .
\end{equation}
Furthermore, from \cite[Lemma 2.2]{[6]}, we have the following compact embedding.

\bl\label{l1}
The space $ X $ is compactly embedded in $ \Ls $ for all $ s \geq 2 $.
\el
Consequently, one can easily verify by standard arguments that $ I $ is well-defined on $ X $, $ I\in C^1(X, \R) $ and
$$
P'(u)(v) = 4 \intR \intR V(|x-y|)u^2(x)u(y)v(y) \ dx dy , \forall \ v \in X .
$$
The derivative of $ P $ can be treated in a similar way as that in \cite[Lemma 2.2]{[6]}. We recall the reader that a nontrivial \textit{weak solution} for \eqref{P} is a function $ u \in X \setminus \{0\} $ satisfying
\begin{align*}
& a \intR \nabla u \nabla v dx + b\left(\intR |\nabla u|^2 dx \right)\intR \nabla u \nabla v dx + \intR Q(x) u v dx \\
& + \mu \intR \intR V(|x-y|)u^2(x)u(y)v(y) \ dx dy = \intR f(u) v dx  , \forall \ v \in X.
\end{align*}
Hence, critical points for $ I $ will be weak solutions for \eqref{P}.

In the sequence we provide two crucial technical lemmas. The first one states when we have boundedness or convergence in $ X $ and the second one shows an important integral convergence. We will prove a general version of theses results, in order to provide a version that can be possible used in some other problems. In this sense, consider a continuous function $ g: \Rdois \RA \R $ satisfying the following condition:

\noindent $
(g) \ g \in C(\R , \R), g(0)=0 \mbox{ and there are constants } K_5\in \R \setminus\{0\}, K_6 > 0 \mbox{ such that } $
$$
K_5 |t| \leq |g(t)| \leq K_6 |t| , \forall \ t \in \R
$$
From condition $ (g) $, we have
\begin{equation}\label{eq1}
\frac{K_5}{2} |t|^2 \leq G(t) \leq \frac{K_6}{2} |t|^2 , \forall \ t \in \R .
\end{equation}

\be\label{ex3}
\textbf{(1)} Clearly, the prototype for $ g $ is given by $ g(t)=t $, for all $ t \in \R $.

\textbf{(2)} We also have more general examples for $g$, such as $ g: \R \RA \R $ given by
$$
g(t) = \left\{ \begin{array}{lll}
t \ , \mbox{ if } t \in [0, 1] \\
t^3 , \mbox{ if } t \in (1, 2] \\
4t , \mbox{ if } t \in (2, \infty)
\end{array} \right. .
$$
One can easily verify that $ g $ satisfies the desired condition.
\ee

\bl\label{l2}
Let $ u\in \Ldois \setminus\{0\} $. Suppose that $ (\un), (\vn) \subset X $ are two sequences satisfying $ \un(x) \RA u(x) $ a.e. in $ \Rdois $ and $ (\vn) $ is bounded in $ \Ldois $. Set
$$
\omega_n = \intR \intR V^+ (|x-y|) G(\un(x))G(\vn(y)) \ dx dy .
$$
Then, if $ \sup\limits_{n \in \N} \omega_n < \infty $, $ (||\vn||_\ast) \subset \R $ is bounded. Moreover, if $ \omega_n \RA 0 $ and $ \vn \RA 0 $ in $ \Ldois $, then $ ||\vn||_\ast \RA 0 $.
\el
\begin{proof}
From Egorov's Theorem, there are $ R \in \N $, $ \delta > 0 $, $ n_0 \in \N $ and $ A \subset B_R $ such that $ A $ is a measurable set with $ |A| > 0 $ and $ \un(x) > \delta $, for all $ n \geq n_0 $. Without loss of generality we can consider $ R > 2 $. Thus, if $ x\in B_R $ and $ y\in B_{2R}^{c} $, we have $ 1+|x-y| \geq \sqrt{1+|y|} $ and $ |x-y|> 2 $. Therefore, for each $ n \geq n_0 $, we have
\begin{align*}
\omega_n & \geq \dfrac{K_{5}^{2}}{4}\ds_{B_{2R}^{c}} \ds_{A} a_1(|x-y|)\ln(1+|x-y|)\un^2(x)\vn^2(y) \ dxdy \\
& \geq \dfrac{K_{5}^{2}a_{1, 0}\delta^2 |A|}{8} \ds_{B_{2R}^{c}} \ln(1+|y|)\vn^2 (y) \ dy \\
& = \dfrac{K_{5}^{2}a_{1, 0}\delta^2 |A|}{8} (||\vn||_{\ast}^{2}-\ln(1+2R)||\vn||_{2}^{2}),
\end{align*}
and the result follows.
\end{proof}

\bl\label{l3}
Let $ (\un)\subset X $ such that $ \un \CF u $ in $ X $. Then,
$$
\lim\limits_{n \RA \infty} \intR \intR V^+ (|x-y|) G(\un(x))g(u(y))(\un(y)-u(y)) \ dxdy = 0.
$$
\el
\begin{proof}
For simplicity, for each $ n \in \N $, set
$$
A_n = \intR \intR V^+ (|x-y|) |G(\un(x))| |g(u(y))| |(\un(y)-u(y))| \ dxdy .
$$
Since $ \un \CF u $ in $ X $, from Lemma \ref{l1}, $ \un \RA u $ in $ \Ls $ for all $ s \geq 2 $. From this, condition $ (g) $ and \eqref{eq1}, we have
\begin{align*}
A_n & \leq ||a_2||_\infty \dfrac{K_{6}^{2}}{2} u^2(x) |u(y)||\un(y) -u(y)| \ dx dy \\
& \leq  \dfrac{K_{6}^{2}||a_2||_\infty}{2} \left[ ||\un||_{\ast}^{2}||u||_{2}^{2}||\un - u ||_{2}^{2} + ||\un||_{2}^{2}\left(\intR \ln(1+|y|) |u(y)||\un(y)-u(y)| \ dy \right)\right].
\end{align*}
To finish the proof one can agree similarly as in \cite[Lemma 2.6]{[6]}.
\end{proof}

\section{Geometry Properties and Convergence Results}

In the present section, we verify that $ I $ has the mountain pass geometry and prove that, up to a subsequence, Cerami sequences are bounded in $ X $. Moreover, we analyse the geometry of the potential $ V $ and the functional $ P $. For the next computations, we need the following useful inequality,
\begin{equation}\label{eq8}
\intR f(u)u \ dx \leq \varepsilon ||u||_{\tau + 1}^{\tau + 1} + K_2||u||_{qr_2}^{q}\left(\intR (e^{r_1 \alpha |u|^2}-1) dx \right)^{\frac{1}{r_1}}.
\end{equation}

\bl\label{l5}
There exists a value $ \rho > 0 $, sufficiently small, such that
$$
m_\beta = \inf \{ I(u) \ ; \ u\in X , ||u|| = \beta \} > 0 \ , \ \forall \beta \in (0, \rho]
$$
and
$$
l_\beta = \inf \{ I'(u)(u) \ ; \ u\in X , ||u|| = \beta \} > 0 \ , \ \forall \beta \in (0, \rho].
$$
\el
\begin{proof}
Let $ \alpha > 4\pi $ and $ u\in X $ such that $ r_1 \alpha ||u||^2 < 4\pi $.

\noindent \textit{Case $a>0$}: From $ (Q) $, equations \eqref{eq7}, \eqref{eq8} and \eqref{eq10}, Lemma \ref{l4} and the embeddings, we have
\begin{align*}
I(u) & \geq \dfrac{a}{2}||\nabla u||_{2}^{2} + \dfrac{b}{4}||\nabla u||_{2}^{4}+\dfrac{Q_0}{2}||u||_{2}^{2}-\dfrac{\mu}{4} K_4||u||_{\frac{8}{3-\lambda}}^{4}-\varepsilon ||u||_{\tau + 1}^{\tau + 1}-K_2 K_\alpha ||u||_{qr_2}^{q} \\
& \geq C_1 ||u||^2 [1-\mu C_2||u||^2 - \varepsilon C_3 ||u||^{\tau-1} - C_4||u||^{q-2}] + \dfrac{b}{4}||\nabla u||_{2}^{4},
\end{align*}
where $ C_1 = \min \left\{ \frac{a}{2}, \frac{Q_0}{2}\right\} >0$ and, similarly,
$$
I'(u)(u) \geq C_5 ||u||^2 [1-\mu C_6||u||^2 - \varepsilon C_7 ||u||^{\tau-1} - C_8||u||^{q-2}] + b||\nabla u||_{2}^{4},
$$
where $ C_5 = \min\{a , Q_0\} $. Hence, for any $ \mu >0 $ and for $ \rho, \varepsilon > 0 $ sufficiently small, the result is valid for this case.

\noindent \textit{Case $a=0$}: Since we are going to take $ \rho >0 $ small, we can assume that $ ||u||^2 < 1 $. Then, $ ||u||_{2}^{2}\geq ||u||_{2}^{4} $. Moreover,
$$
\dfrac{b}{4}||\nabla u||_{2}^{4}+\dfrac{Q_0}{2}||u||_{2}^{4}\geq C_9(||\nabla u||_{2}^{4}+||u||_{2}^{4}) \geq \dfrac{C_9}{4}||u||^4,
$$
where $ C_9=\min\left\{\frac{b}{4} , \frac{Q_0}{2}\right\} $. Consequently,
\begin{equation}\label{eq11}
I(u) \geq ||u||^{4}\left(\dfrac{C_9}{4}-C_{10}\mu -\varepsilon C_{11} ||u||^{\tau - 3}- C_{12}||u||^{q-4}\right) .
\end{equation}
Therefore, for $ \rho, \varepsilon, \mu > 0 $ sufficiently small, we also have the result in this case.
\end{proof}

\bo\label{obs1}
In inequality \eqref{eq11}, we can write the right hand side as
\begin{small}
\begin{equation}\label{eq12}\begin{aligned}
\dfrac{C_9}{4}-C_{10}\mu &-\varepsilon C_{11} ||u||^{\tau - 3}- C_{12}||u||^{q-4}\\
& = \left( \dfrac{C_9}{8}-C_{10}\mu \right) + \left(\dfrac{C_9}{8} -\varepsilon C_{11} ||u||^{\tau - 3}- C_{12}||u||^{q-4}\right).\end{aligned}
\end{equation}
\end{small}
For the first term to be positive we need
$$
\dfrac{C_9}{8C_{10}}> \mu .
$$
To be more precise, we can explicit the constants $ C_9 $ and $ C_{10} $. First of all, from the Gagliardo-Nirenberg inequality, we have
$$
||u||_{\frac{8}{3-\lambda}}^{4}\leq K_{GN}^{\frac{3-\lambda}{2}}||u||^4,
$$
where $ K_{GN}>0 $ is the best constant. Hence, taking into account that
$$
C_{9}= \min\left\{\frac{b}{4} , \frac{Q_0}{2}\right\} \mbox{ \ \ \ and \ \ \ } C_{10}= K_{HLS} K_{GN}^{\frac{3-\lambda}{2}},
$$
we get
$$
\dfrac{\min\left\{\frac{b}{4} , \frac{Q_0}{2}\right\}}{8 K_{HLS}K_{GN}^{\frac{3-\lambda}{2}}} > \mu.
$$
Therefore, considering a sufficiently small $ \mu_0 >0 $ such that the above equation is satisfied and the second term in \eqref{eq12} is positive we have the result of Lemma \ref{l5}, in the \textit{degenerate} case, for all $ \mu \in (0, \mu_0) $.

Finally, we highlight that, in the case which $ a_3 \in \Linf $, we have
$$
\dfrac{\min\left\{\frac{b}{4} , \frac{Q_0}{2}\right\}}{8||a_3||_\infty K_{HLS}K_{GN}^{\frac{3}{2}}} > \mu.
$$
\eo

\bl\label{l6}
Let $ u\in X \setminus\{0\} $ and $ q>4 $. Then,
$$
I(tu)\searrow 0 \mbox{, as \ } t\RA 0 \ , \ \sup\limits_{t>0}I(tu) < \infty \mbox{ \ and \ } I(tu)\RA - \infty \mbox{, as \ } t\RA \infty .
$$
\el
\begin{proof}
Let $ u\in X\setminus\{0\} $, $ q>4 $ and $ t> 0 $. From $ (f_4) $ \eqref{eq13} and \eqref{eq9}, we have
$$
I(tu)\leq \dfrac{a}{2}t^2||\nabla u||_{2}^{2} + \dfrac{b}{4}t^4 ||\nabla u||_{2}^{4}+\dfrac{t^2}{2}||Q||_p ||u||_{2p'}^{2}+\dfrac{\mu K_3}{4}t^4||u||_{X}^{4}-C_q t^q ||u||_{q}^{q} \RA - \infty ,
$$
as $  t\RA \infty $. Now, let $ t>0 $ sufficiently small such that $ r_1 \alpha t^2 ||u||^2 < 4 \pi $. Then, from Lemma \ref{l4} and \eqref{eq7},
$$
\intR F(tu) \ dx \leq \varepsilon t^{\tau + 1}||u||_{\tau + 1}^{\tau + 1} + C_1 t^q ||u||_{qr_2}^{q}\RA 0 \mbox{, as \ } t\RA 0.
$$
Thus, we conclude that $ I(tu)\RA 0 $ as $ t\RA 0 $ and, since $ I\in C^1(X, \R) $, $ \sup\limits_{t>0} I(tu) < \infty $.
\end{proof}

From Lemmas \ref{l5} and \ref{l6}, the value $ c_{mp} $ stated in \eqref{eq16} is well-defined and satisfies $ 0< m_\rho \leq c_{mp} < \infty $. Moreover, since $ I $ has the mountain pass geometry, there exists a Cerami sequence for $ I $ at the level $ c_{mp} $, that is, there exists $ (\un)\subset X $ such that
\begin{equation}\label{eq14}
I(\un)\RA c_{mp} \mbox{ \ \ \ and \ \ \ } ||I'(\un)||_{X'}(1+||\un||_X) \RA 0 \mbox{, as \ } n \RA \infty .
\end{equation}

Before we investigate boundedness and convergence for such sequences, we will study the geometry of $ P $ and $ V $.

\bl\label{l7}
For the potential $ V $ and the functional $ P $, we have the following properties:
\begin{itemize}
\item[(i)] $ V^+ (t)\RA \infty $, as $ t\RA \infty $, and $ V^+ (t)\RA 0 $, as $ t\RA 0 $;

\item[(ii)] $ V^- (t)\RA 0 $, as $ t\RA \infty $;

\item[(iii)] $ V(t)\RA \infty $, as $ t\RA \infty $.

\item[(iv)] There exists a function $ u_0 \in X \setminus\{0\} $ such that $ P(u_0)<0 $.
\end{itemize}
\el
\begin{proof}
\textbf{(i)} From condition $ (V_1) $, for $ t\geq 2 $, we have
$$
0<a_{1, 0}\ln(1+t)<a_1(t)\ln(1+t)\geq V^+(t) \RA 0 \mbox{, as \ } t\RA \infty .
$$
On the other side, also from condition $ (V_1) $,
$$
0 \leq V^+(t) \leq ||a_2||_\infty \ln(1+t)\RA 0, \mbox{ as \ } t\RA 0.
$$
\textbf{(ii)} From condition $ (V_2) $, follows that
$$
0\leq V^-(t)\leq \dfrac{a_3(t)}{t} \leq \left\{\begin{array}{ll}
\dfrac{||a_3||_\infty}{t}, \mbox{ if \ } a_3 \in \Linf \\
\dfrac{1}{t^{1+\lambda}} \ \ \ \ , \mbox{ if \ } a_3(t)=t^{-\lambda}
\end{array}
\right. \RA 0, \mbox{ as \ } t\RA \infty .
$$
\textbf{(iii)} It follows immediately from items (i) and (ii).

\noindent \textbf{(iv)} From condition $ (V_3) $ we can consider an open interval $ (c, d)\subset \I $ in which $ V(t) < 0 $. Take $ x_0 \in (B_{c}^{c}\cap B_{d}) $. Let $ \psi $ be a function such that $ \psi \in C^{\infty}(\Rdois,\R) $ and $ supp \ \psi \subset B_{\frac{|c-d|}{4}}(x_0) $. Then, $ \psi \in X \setminus\{0\} $ and $ P(\psi) < 0 $.
\end{proof}

As an immediately consequence, we have the following corollary.

\bc\label{c1}
The set $ \A = \{ u \in X \ ; \ u \neq 0 , P(u) \leq 0 \} \neq \emptyset $.
\ec

Hence, we are able to find an useful upper bound to the mountain pass level which will make possible to obtain our main results.

\bl\label{l8}
There exists a constant $K_{7}=K_7(a, b, q, Q, p)>0 $ such that $ c_{mp}\leq \dfrac{K_7}{C_{q}^{\frac{2}{q-2}}} $.
\el
\begin{proof}
From the continuous Sobolev embeddings, for $ q>4 $, there exists a constant $ C>0 $ such that $ ||u||\geq C||u||_q $, for all $ u\in \Hum\setminus\{0\} $. Thus, by Corollary \ref{c1}, it makes sense to define
$$
S_q(v)=\dfrac{||v||}{||v||_q} \mbox{ \ \ \ and \ \ \ } S_q = \inf\limits_{v\in \A} S_q(v) \geq \inf\limits_{v\neq 0}S_q(v) > 0 .
$$
Now, from Lemma \ref{l6}, for $ v\in \A $ and $ T>0 $ sufficiently large, $ I(Tv)<0 $. So, we can define a path $ \gamma\in \Gamma $ by $ \gamma(t)=tTv $, for $ t \in [0, 1] $, such that
$$
c_{mp} \leq \max\limits_{0 \leq t \leq 1}I(\gamma(t)) = \max\limits_{0 \leq t \leq 1}I(tTv) \leq \max\limits_{t> 0}I(tv) .
$$
Consequently, from $ (Q) $, $ (f_4) $ and the Gagliardo-Nirenberg inequality, for $ \psi \in \A $, we have
\begin{align*}
c_{mp} & \leq \max\limits_{t>0} \left\{ \left(\dfrac{a+||Q||_p K_{GN}^{\frac{p-1}{p}}}{2}\right)S_q(\psi)^2 t^2 ||\psi||_{q}^{2} - \dfrac{C_q}{2}t^q ||\psi||_{q}^{q}\right\} \\
& + \max\limits_{t>0} \left\{ \dfrac{b}{4}S_q(\psi)^4 t^4 ||\psi||_{q}^{4} - \dfrac{C_q}{2}t^q ||\psi||_{q}^{q} \right\}.
\end{align*}
Considering the auxiliary functions $ h_1 , h_2: \R \RA \R $ given, respectively, by $ h_1(t)= \mathfrak{a}t^2 - \mathfrak{b}t^q $ and $ h_2(t)=\mathfrak{c}t^4 + \mathfrak{d}t^q $, for $ \mathfrak{a}, \mathfrak{b}, \mathfrak{c}, \mathfrak{d} >0 $, we obtain that
$$
c_{mp} \leq \left(2^{\frac{4-q}{q-2}}-\dfrac{2^{\frac{2}{q-2}}}{q}\right)(a+||Q||_p K_{GN}^{\frac{p-1}{p}})^{\frac{q}{q-2}}S_{q}(\psi)^{\frac{2q}{q-2}}\left(\dfrac{1}{qC_q}\right)^{\frac{2}{q-2}}.
$$
Therefore, taking the infimum over all $ \psi \in \A $, we get the desired result.
\end{proof}

Finally, in the last results of this section we verify when Cerami sequences are, up to subsequences, bounded in $ X $. Consider $ (\un)\subset X $ satisfying
\begin{equation}\label{eq15}
\exists \ d >0 \mbox{ \ s.t. \ } I(\un) \leq d , \mbox{ for all } n \in \N \mbox{ and } ||I'(\un)||_{X'}(1+||\un||_X) \RA 0 , \mbox{ as } n \RA \infty .
\end{equation}

\bl\label{l9}
Let $ (\un)\subset X $ be bounded in $ \Hum $ such that
\begin{equation*}
\liminf\limits_{n\RA \infty} \sup\limits_{y\in \mathbb{Z}^2} \ds_{B_2(x)}\un^2(x) dx > 0 .
\end{equation*}
Then, there exists $ u\in \Hum \setminus\{0\} $ and $ (\yn)\subset \mathbb{Z}^2 $ such that, up to a subsequence, $ \yn \ast \un = \until \CF u\in \Hum $. Particularly, $ u\neq 0 $ in $ \Ldois $.
\el

\bl\label{l11}
Let $ (\un)\subset X $ be a sequence satisfying (\ref{eq15}), bounded in $ \Hum $ and such that $ ||\nabla \un ||_2 < 2\sqrt{\frac{\pi}{r_1 \alpha}} $, for all $ n \in \N $, and
$$
\liminf_{n\RA \infty} \sup\limits_{y\in \mathbb{Z}^2} \ds_{B_2(y)}\un^2(x) dx > 0 .
$$
Then, up to a subsequence, $ (\until) $ is bounded in $ X $.
\el
\begin{proof}
The proof follows from Lemmas \ref{l1}, \ref{l2} and \ref{l9}, equations \eqref{eq7} and \eqref{eq10} and the facts that, $ P_1 $ is invariant under $ \mathbb{Z}^2$-translations and that, for all $ n \in \N $,
$$
\dfrac{\mu}{4}P_1(\un)  = I(\un) - \dfrac{a}{2}||\nabla \un||_{2}^{2} - \dfrac{b}{4}||\nabla \un ||_{2}^{4} - \dfrac{1}{2}\intR Q(x) \un^2 (x) dx + \dfrac{\mu}{4}P_2(\un) + \int F(\un) dx,
$$
as desired.
\end{proof}

We highlight that the next technical lemma is the key in obtaining multiplicity of solutions for problem \eqref{P},  since it makes possible to verify the validity of $(PS)$ condition at some suitable levels.

\bc\label{l12}
Let $ (\un)\subset X $ under the hypotheses given in Lemma \ref{l11}. Then, up to a subsequence, $ (\un) $ is bounded in $ X $.
\ec
\begin{proof}
To begin with, from Lemma \ref{l11}, passing to a subsequence if necessary, there exists $ (\yn)\subset \mathbb{Z}^2 $ such that $ \until \CF u $ in $ X $, with $ u \neq 0 $ in $\Ldois$, $ \until(x) \RA u(x) $ pointwise a.e. in $ \Rdois $ and, from Lemma \ref{l1}, $ \until \RA u $ in $ \Ls $, for all $ s \geq 2 $.

Moreover, one can see that there are $ R_1, C_1 >0 $ and $ n_1 \in \N $ such that $ ||\un||_{p, B_{R_1}}^{p} \geq C_1 > 0 $, for all $ n \geq n_1 $. From this we can conclude that $ (\yn) $ is bounded in $ \mathbb{Z}^2 $ and, using that
$$
||\un||_{\ast}^{2} = \intR \ln(1+|x- \yn|) \until^2(x) dx \leq ||\until||_{\ast}^{2}+\ln(1+|\yn|)||\until||_{2}^{2} , \forall \ n \in \N,
$$
and that $ (\un) $ is already bounded in $ \Hum $, the result follows.
\end{proof}

\bl\label{l10}
Assume $ q>4 $ and $ \alpha > 4\pi $, fixed. Let $ (\un)\subset X $ be a sequence satisfying (\ref{eq15}), $ ||\nabla \un ||_2 < 2\sqrt{\frac{\pi}{r_1 \alpha}} $ and that does not verify $ ||\un||\RA 0 $ and $ I(\un) \RA 0 $. Then,
$$
\liminf_{n\RA \infty} \sup\limits_{y\in \mathbb{Z}^2} \ds_{B_2(y)}\un^2(x) dx > 0 .
$$
\el
\begin{proof}
The proof is done by contradiction, applying the Lion's Lemma, and using \eqref{eq10}, \eqref{eq8}, Moser-Trudinger inequality and that $ I'(\un)(\un)\RA 0 $, as $ n \RA \infty $.
\end{proof}

\section{The \textit{nondegenerate} case ($a>0$)}

This section is devoted to prove Theorems \ref{t1} and \ref{t2}. Since we are going to handle the \textit{nondegerate} case, throughout this section we will assume $ a>0 $ and $ b \geq 0 $. Our strategy consists in proving boundedness of Cerami sequences in $\Hum$, guaranteeing that it is possible to apply Moser-Trudinger inequality for such sequences and, under what conditions, $ I $ has nontrivial critical points in $ X $. To finish this section, we verify that \eqref{P} has infinitely many solutions.

\bl\label{l13}
Suppose that $ a>0 $ and $ b\geq 0 $. Let $ (\un)\subset X $ a sequence satisfying \eqref{eq15}. Then, $ (\un) $ is bounded in $ \Hum $.
\el
\begin{proof}
From condition $ (f_3) $ and \eqref{eq15}, we have
\begin{align*}
d+ o(1) & \geq I(\un) - \dfrac{1}{4}I'(\un)(\un) = \dfrac{a}{4}||\nabla \un||_{2}^{2}+\dfrac{1}{4}\intR Q(x)\un^2 dx +\intR \left[\dfrac{f(\un)\un}{4}-F(\un) \right] dx \\
& \geq \dfrac{\min\{a, Q_0\}}{4}||\un||^2 , \forall \ n \in \N.
\end{align*}
Hence, for all $ n \in \N $,
$$
\left(\dfrac{4d}{\min\{a, Q_0\}}\right)^{\frac{1}{2}}+o(1) \geq ||\un|| , \forall \ n\in \N,
$$
and the result follows.
\end{proof}

\bc\label{c2}
Let $ (\un)\subset X $ be a sequence satisfying \eqref{eq15}, with $ d\in (0, c_{mp}] $, or being a Cerami sequence in level $ c_{mp} $. Then, up to a subsequence, there exists a constant $ K_8=K_8(a, b, q, Q, p)>0 $ such that $ ||\un||\leq \dfrac{K_8}{C_{q}^{\frac{1}{q-2}}} $, for all $ n \in \N $.
\ec
\begin{proof}
The proof follows directly from Lemmas \ref{l8} and \ref{l13} and $ \limsup $ properties.
\end{proof}

\bp\label{p1}
Suppose $ q>4 $ and $ C_q > 0 $ sufficiently large. Let $ (\un)\subset X $ a sequence satisfying \eqref{eq15}, with $ d\in (0, c_{mp}] $, or being a Cerami sequence in level $ c_{mp} $. Then, passing to a subsequence, if necessary, only one between the following alternatives hold:
\begin{itemize}
\item[(a)] $ ||\un||\RA 0 $ and $ I(\un)\RA 0 $.

\item[(b)] There exists a function $ u\in X\setminus\{0\} $ such that $ \un \RA u $ in $ X $ and $ u $ is a critical point to $ I $ in $ X $.
\end{itemize}
\ep
\begin{proof}
Let us suppose that item (a) does not hold. Then, from Lemmas \ref{l13}, \ref{l9}, \ref{l10}, \ref{l11} and Corollary \ref{l12}, passing to a subsequence if necessary, $ \un \CF u $ in $ X $, for $ u\in X\setminus\{0\} $. Moreover, from Lemma \ref{l1}, $ \un \RA u $ in $ \Ls $, for all $ s \geq 2 $.

Now, from Corollary \ref{c2}, up to a subsequence, we can assume that $ r_1 \alpha ||\un||^2 < 4\pi $, for all $ n \in \N $ and $ C_q > 0 $ sufficiently large. Thus, from \eqref{eq15}, \eqref{eq13}, Lemma \ref{l4}   and (HLS), we have the following
main properties

\noindent \textbf{(i)} $ |I'(\un)(\un -u)|\leq ||I'(\un)||_{X'}||\un - u||_X \RA 0 $, as $ n\RA \infty $;

\noindent \textbf{(ii)} $ P_2'(\un)(\un - u) \RA 0 $, $ \intR Q(x) \un^2 dx \RA 0 $ and $ \intR f(\un)\un dx \RA 0 $, as $ n \RA \infty $.

Moreover, from $ \un \CF u $ in $ \Hum $, the weakly sequentially lower semicontinuity of $ ||\cdot||_2 $ and $ \liminf $ properties, passing to a subsequence if necessary, we have
$$
\langle \nabla \un , \nabla (\un -u) \rangle = ||\nabla \un||_{2}^{2}-||\nabla u||_{2}^{2}+o(1).
$$
Consequently, from Lemma \ref{l3} with $ g(t)=t $, (i) and (ii), follows that
\begin{align*}
o(1) & = I'(\un)(\un -u) \\
& \geq a (||\nabla \un||_{2}^{2}-||\nabla u||_{2}^{2}) + b||\nabla \un||_{2}^{2}(||\nabla \un||_{2}^{2}-||\nabla u||_{2}^{2}+o(1)) P_{1}'(\un)(\un-u) +o(1) \\
& = a (||\nabla \un||_{2}^{2}-||\nabla u||_{2}^{2}) + \intR \intR V^+(|x-y|)\un^2(x)(\un-u)^2(y) dx dy \\
& + \intR \intR V^+(|x-y|)\un^2(x)u(y)(\un(y)-u(y)) dxdy + o(1) \\
& \geq a (||\nabla \un||_{2}^{2}-||\nabla u||_{2}^{2}) + o(1) \geq o(1).
\end{align*}
Hence, we obtain that $ ||\nabla \un||_{2}^{2} - ||\nabla u||_{2}^{2} \RA 0 $ and, since $ \un \RA u $ in $ \Ldois $, $ \un \RA u $ in $ \Hum $. Moreover, returning to the above inequality we conclude that
$$
\intR \intR V^+(|x-y|)\un^2(x)(\un-u)^2(y) dx dy \RA 0
$$
and, from Lemma \ref{l2}, $ ||\un - u||_\ast \RA 0 $, which implies that $ \un \RA u $ in $ X $. Finally, for $ v \in X $,
$$
|I'(u)(v)| \leq |I'(u)(v)-I'(\un)(v)|+||I'(\un)||_{X'}||v||\RA 0, \mbox{ as } n \RA \infty .
$$
Therefore, $ u $ is a nontrivial critical point for $ I $ in $ X $.
\end{proof}

\begin{proof}[Proof of Theorem \ref{t1}]
Item (a) follows immediately from \eqref{eq14}, Lemma \ref{l5} and Proposition \ref{p1}. Let us prove item (b). From item (a), $ \K \neq \emptyset $. Let $ (\un)\subset \K $ such that $ I(\un)\RA c_g $.

Observe that $ c_g \in [-\infty, c_{mp}] $. If $ c_g=c_{mp} $ nothing remains to be proved. Assume that $ c_g < c_{mp} $. Thus, combined with the definition of $ \K $, we have that $ (\un) $ satisfies \eqref{eq15} with $ d=c_{mp} $. Hence, from Lemma \ref{l5} and Proposition \ref{p1}, there exists $ u\in X\setminus\{0\} $ such that $ \un \RA u $ in $ X $ and $ u $ is a critical point for $ I $. Moreover, we have $ I(u)=c_g $ which implies, particularly, that $ c_g>-\infty $.
\end{proof}

In order to prove our second main result, let $ k\in \N $, arbitrary but fixed, and $ Z\subset X $ a subspace with $ \dim Z =k $ and norm $ ||\cdot||_Z $. Our goal is to apply a symmetric version of the mountain pass theorem, due to Ambrosetti and  Rabinowitz \cite{[ambrosseti]}(see also \cite{[bartolo], [silva]}).

\bt\label{t4}
(\cite[Theorem 4.1]{[albuquerque]}) Let $ E= E_1 \oplus E_2 $, where $ E $ is a real Banach space and $ E_1 $ is finite dimensional. Suppose that $ J \in C^{1}(E, \R) $ is even, $ J(0)=0 $, and that it verifies
\begin{itemize}
\item[$ (J_1) $] there exists $ \tau , r > 0 $ such that $ J(u) \geq \tau $ if $ ||u||_E = r $, $ u\in E_2 $,

\item[$ (J_2) $] there exists a finite-dimensional subspace $ \F \subset E $, with $ \dim E_1 < \dim \F $, and a constant $ \B > 0 $ such that $ \max\limits_{u\in \F} J(u) \leq \B $,

\item[$ (J_3) $] $ J $ satisfies the $ (PS)_c $ condition for all $ c\in (0, \B) $.
\end{itemize}
Then, $ J $ possess at least $ \dim \F - \dim E_1 $ pairs of nontrivial critical points.
\et

In the sequence we need to verify the conditions of Theorem \ref{t4}. First of all, one should observe that, under conditions $ (f_1 ')-(f_4) $, $ (Q) $, $ (M) $ and $ (V_1)-(V_3) $ we already have that $ I \in C^{1}(X, \R) $, is even, $ I(0)=0 $ and, from Lemma \ref{l5}, $ I $ verifies $ (J_1) $. So, it remains to prove that $ I $ also verifies $ (J_2) $ and $ (J_3) $.

\bl\label{l14}
Let $ q>4 $. Then, there exists $ R>0 $ such that $ I(u) \leq 0 $ for all $ u\in X $ verifying $ ||u||_Z \geq R $.
\el
\begin{proof}
Since $ \dim Z < \infty $, all norms are equivalent. Thus, from condition $ (f_4) $ and \eqref{eq9}, we have
$$
I(u) \leq C_1 ||u||_{Z}^{2}+C_2 ||u||_{Z}^{4} - C_3 ||u||_{Z}^{q} \RA - \infty , \mbox{ as } ||u||_Z \RA \infty .
$$
\end{proof}

\bl\label{l15}
Let $ q>4 $. Then, there exists $ \eta > 0 $, sufficiently small, such that $ \max\limits_{u\in Z}I(u) \leq \eta $ and $ r_1 \alpha \frac{\eta}{\min\{a, Q_0\}} < \pi $.
\el
\begin{proof}
Let $ u\in Z \setminus \{0\} $. Thus, from $ \dim Z < \infty $, condition $ (f_4) $ and \eqref{eq9}, there are constants constants $ C_1, C_2, C_3 >0 $, depending on $ a, b, q $ and $ Q $ such that
$$
I(u) \leq C_1 ||u||_{Z}^{2}+ C_2 ||u||_{Z}^{4} - C_q C_3 ||u||_{Z}^{q}.
$$
Arguing in a similar way as in Lemma \ref{l8}, one can find a constant $ C_4>0 $ satisfying
$$
I(u)\leq \dfrac{C_4}{C_{q}^{\beta}}, \mbox{ for some exponent } \beta = \beta(q)>1.
$$
Consequently,
$$
\max\limits_{u\in Z} I(u) \leq \dfrac{C_4}{C_{q}^{\beta}}
$$
and, taking $ C_q > 0 $ sufficiently large we find a value $ \eta >0 $ sufficiently small as desired.
\end{proof}

In the next proposition we guarantee that $ I $ satisfies the $ (PS)_d $ condition for all $ d \in (0, \eta) $. One can observe that the proof can be done in a very similar way as that of Proposition \ref{p1}, so we will omit it here. We highlight that the validity of following lemma is possible only in virtue of Lemma \ref{l12}.

\bl\label{l16}
The functional $ I $ satisfies condition $ (PSC)_d $ for all $ d\in (0, \eta) $.
\el

\begin{proof}[Proof of Theorem \ref{t2}]
From Lemmas \ref{l5}, \ref{l15} and \ref{l16} and an immediate application of Theorem \ref{t4}, with $ E=X $, $ E_1 = \{0\} $, $ \F = Z $, $ J=I $, $ \tau = m_\rho $, $ r= \rho $ and $ \B = \eta $, we get that $ I $ possess at least $ k $ nontrivial critical points. Therefore, as we can make $ k $ as large as we want, we conclude that (\ref{P}) has infinitely many solutions.
\end{proof}

\section{The \textit{degenerate} case ($a=0$)}

In this section we investigate the existence of solutions for \eqref{P} in the \textit{degenerate} case. So, we assume $ a=0 $ and $ b>0 $ throughout it. Since we have the same multiplying constant $ \frac{1}{4} $ in both terms, that one depending on $ ||\nabla \cdot||_2 $ and that on with $ V $, we need a different approach than was used in Section 4. The technique is based in Lemmas \ref{l17} and \ref{l19}, which are inspired in similar results of \cite{[6]}.

\bl\label{l17}
Let $ (\un)\subset X $ a sequence satisfying \eqref{eq15} and $ (t_n)\subset \left. \left(0, \left(\frac{\theta - 4}{\theta}\right)^{\frac{1}{4}}\right]\right.$. Then, $ I(t_n \un) \leq I(\un) $, for all $ n \in \N $.
\el
\begin{proof}
Observe that
\begin{equation}\label{eq17}
I(t_n \un) = \dfrac{b}{4}t_{n}^{4}||\nabla \un ||_{2}^{4}+\dfrac{t_{n}^{2}}{2}\intR Q(x)\un^2 dx + \dfrac{\mu}{4}t_{n}^{4}P(\un)-\intR F(t_n \un) dx
\end{equation}
and
\begin{equation}\label{eq18}
\mu P'(\un)(\un) = I'(\un)(\un)-b||\nabla \un||_{2}^{4}-\intR Q(x) \un^2 dx + \intR f(\un) \un dx .
\end{equation}
Thus, combining \eqref{eq17} and \eqref{eq18},
\begin{small}
\begin{align*}
I(t_n \un) - I(\un) & = \dfrac{b}{4}(t_{n}^{4}-1)||\nabla \un ||_{2}^{4}+\dfrac{t_{n}^{2}-1}{2}\intR Q(x)\un^2 dx + \dfrac{\mu}{4}(t_{n}^{4}-1)P(\un)\\
&\hspace{3cm}+\intR [F(\un)- F(t_n \un)] dx \\
& = \dfrac{1}{2}\left(t_{n}^{2}-\dfrac{t_{n}^{4}}{2}-\dfrac{1}{2}\right)\intR Q(x) \un^2 dx + \intR \left[F(\un)-F(t_n \un) + \dfrac{t_{n}^{4}-1}{4}f(\un)\un \right] dx \\
& \leq 0,
\end{align*}
\end{small}
since $ t_{n}^{2}-\frac{t_{n}^{4}}{2}-\frac{1}{2}\leq 0 $, for all $ n \in \N $, and
$$
F(\un)+\dfrac{t_{n}^{4}-1}{4}f(\un)\un \leq \left(\dfrac{1}{\theta}+\dfrac{t_{n}^{4}-1}{4}\right)f(\un)\un \leq 0, \forall \ n \in \N.
$$
Therefore, $ I(t_n \un) \leq I(\un) $, for all $ n \in \N $.
\end{proof}

\bl\label{l19}
Let $ (\un)\subset X $ satisfying \eqref{eq15}. Then, $ (\un) $ is bounded in $ \Hum $.
\el
\begin{proof}
Suppose by contradiction that $ ||\un||\RA \infty $. For a fixed $ \alpha > 4\pi $, define $ \vn = \sqrt{\frac{\pi}{r_1 \alpha}}\frac{\un}{||\un||} $, for each $ n \in \N $. Thus, $ ||\vn||=\sqrt{\frac{\pi}{r_1 \alpha}} $, for all $ n \in \N $. Consequently, $ (\vn) $ is bounded in $ \Hum $.

\noindent \textbf{Claim:} $ \liminf\limits_{n \RA \infty}\sup\limits_{y \in \mathbb{Z}^2}\ds_{B_2(y)}\vn^2(x) dx >0 $.

Otherwise, by Lion's Lemma, $ \vn \RA 0 $ in $ \Ls $, for all $ s \in (2, \infty) $. Thus, from \eqref{eq10}, \eqref{eq8}, \eqref{eq15} and Lemma \ref{l4}, we have
\begin{align*}
0 \leq \mu P_1(\vn)+b||\nabla \vn||_{2}^{4}+\dfrac{1}{2}Q(x)\vn^2 dx & = I'(\vn)(\vn)+\mu P_2(\vn)+\intR f(\vn)\vn dx \\
& \leq \mu K_5 ||\vn||_{\frac{8}{3-\lambda}}^{4}+\varepsilon ||\vn||_{\tau + 1}^{\tau + 1} + C_1 ||\vn||_{qr_2}^{q} \RA 0,
\end{align*}
as $ n \RA \infty $. Consequently, $ P_1(\vn) \RA 0 $, $ ||\nabla \vn||_2 \RA 0 $ and, by $ (Q) $, $ ||\vn||_2 \RA 0 $. Hence, $ ||\vn||\RA 0 $, which is a contradiction.

Therefore, from Lemma \ref{l12}, up to a subsequence, $ \vn \CF v $ in $ X $, for $ v \in X \setminus\{0\} $. We can assume, without loss of generality, that $ \vn(x)\RA v(x) $ a.e. in $ \Rdois $. Moreover, from the continuous Sobolev embeddings, $ \vn \CF v  $ in $ L^q(\Rdois) $. Thus, from weakly sequentially lower semicontinuity, boundedness in $ X $, \eqref{eq13}, the Gagliardo-Nirenberg inequality and condition $ (f_4) $, there exist $ n_0\in \N $ such that, for $ t>0 $,
$$
I(t\vn) \leq t^4 \dfrac{b}{4}\left(\dfrac{\pi}{r_1 \alpha}\right)^2 + \dfrac{t^2}{2}\dfrac{\pi}{r_1 \alpha}||Q||_p K_{GN}^{\frac{1}{p'}}+ t^4 C_2 -C_q t^q ||v||_{q}^{q}, \forall \ n \geq n_0.
$$
So, if we choose $ t_0 > 0 $ sufficiently large, $ I(t_0 \vn) \leq -1 $, for all $ n\geq n_0 $. But, by other hand,
$$
t_0 \sqrt{\dfrac{\pi}{r_1 \alpha}} \dfrac{1}{||\un||}\RA 0, \mbox{ as } n \RA \infty,
$$
which contradicts Lemma \ref{l6}. Therefore, $ (\un) $ is bounded in $ \Hum $.
\end{proof}

Let $ (\un)\subset X $ be the sequence given in equation \eqref{eq14}. From Lemma \ref{l19}, there exists a constant $ K_{mp}>0 $ such that, passing to a subsequence ie necessary, $ ||\un||\leq K_{mp} $, for all $ n \in \N $. Although we already have a bound for this Cerami sequence, we still need a sufficiently small bound for $ ||\nabla \un||_2 $ in order to apply Lemma \ref{l4}.

\bl\label{l20}
Let $ (\un)\subset X $ a sequence satisfying \eqref{eq14}. Then, there are $ \mu_{mp} > 0 $ sufficiently small and a constant $ K_9 > 0 $ such that, up to a subsequence, $ ||\nabla \un||_2 \leq K_9 c_{mp}^{\frac{1}{4}} $, for all $ n \in \N $ and $ \mu \in (0, \mu_{mp}) $.
\el
\begin{proof}
From \eqref{eq10} and the Gagliardo-Nirenberg inequality, we have
$$
P_2(\un) \leq K_5 ||\un||_{\frac{8}{3-\lambda}}^{4}\leq K_5 K_{GN}^{\frac{3-\lambda}{2}}||\un||^4 \leq K_5 K_{GN}^{\frac{3-\lambda}{2}} K_{mp}^{4}, \forall \ n \in \N.
$$
Consequently,
\begin{align*}
c_{mp}+o(1) & \geq I(\un) - \dfrac{1}{8}I'(\un)(\un) \\
& \geq \dfrac{b}{8}||\nabla \un||_{2}^{4}+ \dfrac{3}{8}Q_0 ||\un||_{2}^{2}-\mu \dfrac{K_5 K_{GN}^{\frac{3-\lambda}{2}} K_{mp}^{4}}{4}+ \intR \left(\dfrac{1}{8}f(\un)\un - F(\un) \right) dx \\
& \geq \dfrac{b}{8}||\nabla \un||_{2}^{4}-\mu \dfrac{K_5 K_{GN}^{\frac{3-\lambda}{2}} K_{mp}^{4}}{4}.
\end{align*}
Hence, considering $ \mu_{mp}>0 $ sufficiently small, we have
$$
c_{mp}+ o(1)\geq \mathfrak{e} b ||\nabla \un||_{2}^{4},
$$
for all $ n \in \N $, $ \mu \in (0, \mu_{mp}) $ and a value $ \mathfrak{e}\in \left(0, \frac{1}{8}\right)  $, depending on $ \mu_{mp} $.

Therefore, the result follows for $ K_9 = \left(\frac{1}{\mathfrak{e} b}\right)^{\frac{1}{4}} > 0 $.
\end{proof}

\begin{proof}[Proof of Theorem \ref{t3} - (a)]
The proof follows from Lemmas \ref{l3}, \ref{l4}, \ref{l5}, \ref{l9}, \ref{l10}, \ref{l11}, \ref{l19} and \ref{l20} and Corollary \ref{l12}, arguing in a very similar way as Proposition \ref{p1} and Theorem \ref{t1}.
\end{proof}

Since we already have item (a), we can consider the set $ \K = \{ v \in X \setminus\{0\} \ ; \ I'(v) = 0\} $ that is not empty. So, $ c_{g}\in [-\infty, c_{mp}] $ and there exists a sequence $ (\un)\subset \K $ such that
$
I(\un) \RA c_g
$. We will assume $ c_g < c_{mp} $ once if the equality holds, nothing remains to be proved. Also, by the definition of $ \K $, one can see that $ (\un) $ satisfies
$$
||I'(\un)||_{X'}(1+||\un||_X)\RA 0, \mbox{ as } n \RA \infty.
$$
Thus, from Lemma \ref{l19}, there exists a constant $ K_g > 0 $ such that $ ||\un||\leq K_g $, for all $ n \in \N $.

\bl\label{l21}
Let $ (\un)\subset X $ be the minimizing sequence for $ c_g $. Then, there are $ \mu_{g} > 0 $ sufficiently small and a constant $ K_{10} > 0 $ such that, up to a subsequence, $ ||\nabla \un||_2 \leq K_{10} c_{mp}^{\frac{1}{4}} $, for all $ n \in \N $ and $ \mu \in (0, \mu_{g}) $.
\el
\begin{proof}
Similarly as done in Lemma \ref{l20}, since $ c_g < c_{mp} $, up to a subsequence, we get
$$
c_{mp}+o(1) \geq \dfrac{b}{8}||\nabla \un||_{2}^{4}-\mu \dfrac{K_5 K_{GN}^{\frac{3-\lambda}{2}} K_{g}^{4}}{4}.
$$
Thus, considering $ \mu_{g}>0 $ sufficiently small, we have
$$
c_{mp}+ o(1)\geq \mathfrak{r} b ||\nabla \un||_{2}^{4},
$$
for all $ n \in \N $, $ \mu \in (0, \mu_{g}) $ and a value $ \mathfrak{r}\in \left(0, \frac{1}{8}\right)  $, depending on $ \mu_{g} $.

Therefore, the result follows for $ K_{10} = \left(\frac{1}{\mathfrak{r} b}\right)^{\frac{1}{4}} > 0 $.
\end{proof}

\bo\label{obs2}
One should observe that, since the nonemptiness of $ \K $ depends on the existence of a solution at the mountain pass level, the value $ \mu_g $ must also satisfy $ \mu_g \leq \mu_{mp} $.
\eo

\begin{proof}[Proof of Theorem \ref{t3}-(b)]
The proof follows from Lemmas \ref{l3}, \ref{l4}, \ref{l5}, \ref{l9}, \ref{l10}, \ref{l11}, \ref{l19} and \ref{l21} and Corollary \ref{l12}, arguing in a very similar way as Proposition \ref{p1} and Theorem \ref{t1}.
\end{proof}

\noindent \textbf{Acknowledgements:} The first author was supported by  Coordination of Superior Level Staff Improvement-(CAPES)-Finance Code 001 and  S\~ao Paulo Research Foundation-(FAPESP), grant $\sharp $ 2019/22531-4, while the second  author was supported by  National Council for Scientific and Technological Development-(CNPq),  grant $\sharp $ 307061/2018-3 and FAPESP  grant $\sharp $ 2019/24901-3. The third author 
is a member of the {\em Gruppo Nazionale per
l'Analisi Ma\-te\-ma\-ti\-ca, la Probabilit\`a e le loro Applicazioni} (GNAMPA)
of the {\em Istituto Nazionale di Alta Matematica} (INdAM) and was partly supported by
the {\em Fondo Ricerca di Base di Ateneo -- Eser\-ci\-zio 2017--2019} of the University of Perugia, named {\em  PDEs and Nonlinear Analysis}.


\begin{thebibliography}{2}

\bibitem{[albuquerque]} Albuquerque, F. S. B. (2014) Nonlinear Schrodinger elliptic systems involving exponential critical growth in $\mathbb{R}^2$, \textit{Electronic Journal of Differential Equations}. Vol. 2014, n. 59, pp. 1-12.

\bibitem{[ambrosseti]} Ambrosetti, A. and Rabinowitz P.H. (1973) Dual variational methods in critical point theory and applications, \textit{Journal of Functional Analysis}. 14, 349–381.

\bibitem{[bartolo]} Bartolo, P., Benci, V., Fortunato, D. (1983) Abstract critical point theorems and applications to some nonlinear problems with “strong” resonance at infinity, \textit{Nonlinear Analysis: Theory, Methods \& Applications}. 7, 981–1012.

\bibitem{[boer]} Böer, E. de S. and Miyagaki, O. H. (2021) Existence and multiplicity of solutions for the fractional $p$-Laplacian Choquard logarithmic equation involving a nonlinearity with exponential critical and subcritical growth, \textit{J. Math. Phys.} 62, 051507.

\bibitem{[boer2]} Böer, E. de S. and Miyagaki, O.H. (2021) $ (p, N)$-Choquard logarithmic equation involving a nonlinearity with exponential critical growth: existence and multiplicity, \textit{ArXiv:2105.11442 [Math]}, submitted.

\bibitem{[5]} Cao, D. M. (1992) Nontrivial solution of semilinear elliptic equations with critical exponent in $\mathbb{R}^2$. \textit{Communications in Partial Differential Equations}, 17, 407–435.

\bibitem{[6]} Cingolani, S. and Weth, T. (2016) On the planar Schrödinger–Poisson system. \textit{Annales de l’Institut Henri Poincare (C) Non Linear Analysis}, 33, 169–197.

\bibitem{[10]} Du, M. and Weth, T. (2017) Ground states and high energy solutions of the planar Schrödinger–Poisson system. \textit{Nonlinearity}, 30, 3492–3515.

\bibitem{[11]} Fröhlich, H. (1937) Theory of electrical breakdown in Ionic crystals. \textit{Proceedings of the Royal Society of London. Series A, Mathematical and Physical Sciences}, 160, 230–241.

\bibitem{[12]} Fröhlich, H. (1954) Electrons in lattice fields. \textit{Advances in Physics}, 3, 325–361.
35

\bibitem{[jin]} Jin, J. and Wu, X. (2010) Infinitely many radial solutions for Kirchhoff-type problems in RN, \textit{Journal of Mathematical Analysis and Applications}. 369, 564–574.

\bibitem{[Lam]} Lam, N. and Lu, G. (2014) Elliptic equations and systems with subcritical and critical exponential growth without the Ambrosetti-Rabinowitz condition, \textit{J. Geom. Anal.}, 24,118-143.

\bibitem{[pucci]} Liang, S., Pucci, P. and Zhang, B. (2020) Multiple solutions for critical Choquard-Kirchhoff type equations, \textit{Advances in Nonlinear Analysis}. 10, 400–419.


\bibitem{[liang]} Liang, S. and Zhang, J. (2015) Existence of solutions for Kirchhoff type problems with critical nonlinearity in $\mathbb{R}^N$, \textit{Z. Angew. Math. Phys.} 66, 47–562.

\bibitem{[15]} Lieb, E. H. (1983) Sharp constants in the Hardy-Littlewood-Sobolev and related inequalities. \textit{The Annals of Mathematics}, 118, 349.

\bibitem{[16]} Lions, P.-L. (1987) Solutions of Hartree-Fock equations for Coulomb systems. \textit{Communications in Mathematical Physics}, 109, 33–97.

\bibitem{[olimpio]} Miyagaki, O.H. and Pucci, P. (2019) Nonlocal Kirchhoff problems with Trudinger–Moser critical nonlinearities, \textit{Nonlinear Differ. Equ. Appl.}, 26:27, 1-26.

\bibitem{[17]} Moser, J. (1971) A Sharp form of an inequality by N. Trudinger. \textit{Indiana University
Mathematics Journal}, 20, 1077–1092.

\bibitem{[18]} Penrose, R. (1996) On gravity’s role in quantum state reduction. \textit{General Relativity and Gravitation}, 28, 581–600.

\bibitem{[pucci2]} Pucci, P., Xiang, M. and Zhang, B. (2019) Existence results for Schrödinger–Choquard–Kirchhoff equations involving the fractional p-Laplacian, \textit{Advances in Calculus of Variations}, 12, 253–275.

\bibitem{[ruf]} Ruf, B. and Sani, F. (2013) Ground states for elliptic equations in $\mathbb{R}^{2}$ with exponential critical growth. Magnanini, R., Sakaguchi, S., and Alvino, A. (eds.), \textit{Geometric properties for parabolic and elliptic PDE’s}, vol. 2, pp. 251–267, Springer, Milan.

\bibitem{[silva]} Silva, E. A. B. (1988) \textit{Critical point theorems and applications to differential equations}, PhD. Thesis, University of Wisconsin-Madison.

\bibitem{[22]} Wilson, A. J. C. (1955) Untersuchungen über die Elektronentheorie der Kristalle by S. I.
Pekar. \textit{Acta Crystallographica}, 8, 70–70.

\end{thebibliography}
\end{document}